\documentclass[12pt]{amsart}

\usepackage{amsmath,xspace,amssymb,mathrsfs}
\usepackage{color}

\input xy
\xyoption{all}
\xyoption{2cell}
\UseAllTwocells
\CompileMatrices

\newcommand{\Spec}{\operatorname{Spec}}
\renewcommand{\phi}{\varphi}

\newcommand{\rn}{\operatorname{rank}}

\newcommand{\MA}{\operatorname{Max}}

\newcommand{\Min}{\operatorname{Min}}
\newcommand{\Ann}{\operatorname{Ann}}

\newcommand{\Supp}{\operatorname{Supp}}
\newcommand{\sgn}{\operatorname{sgn}}

\newtheorem{proposition}{Proposition}[section]
\newtheorem{lemma}[proposition]{Lemma}

\newtheorem{corollary}[proposition]{Corollary}
\newtheorem{theorem}[proposition]{Theorem}

\theoremstyle{definition}

\newtheorem{remark}[proposition]{Remark}

\begin{document}

\title[Upper topology and projective modules]{The upper topology and its relation with the projective modules}

\author[A. Tarizadeh]{Abolfazl Tarizadeh}
\address{Department of Mathematics, Faculty of Basic Sciences, University of Maragheh \\
P. O. Box 55136-553, Maragheh, Iran.
 }
\email{ebulfez1978@gmail.com}

\date{}

\begin{abstract} In this paper, first we obtain some new and interesting results on projective modules and on the upper topology of an ordinal number. Then it is shown that the rank map of a locally of finite type projective module is continuous with respect to the upper topology (by contract, it is well known this map is not necessarily continuous with respect to the discrete topology). It is also proved that a finitely generated flat module is projective if and only if its rank map is continuous with respect to the upper topology.
\end{abstract}

\maketitle

\section{Introduction}

This paper grew out of an attempt to understand when the rank map of an infinitely generated projective module is continuous. We use f.g. in place of ``finitely generated". It is well known that the rank map of a f.g. projective $R-$module is continuous whenever $\Spec(R)$ is equipped with the Zariski topology and the set of natural numbers with the ``discrete" topology. In Proposition \ref{proposition 2}, this fact is also re-proved by a new method, see e.g. \cite[Corollary 2.2.2]{Weibel} for another proof of it. Note that the ``finitely generated'' assumption of  module is a crucial point in the continuity of the rank map. If we drop this hypothesis then the rank map is not necessarily continuous. For example, there are locally of finite type projective modules whose rank maps are not necessarily continuous with respect to the discrete topology, see e.g. \cite[Ex. 2.15]{Weibel}. Therefore it is natural to ask, under which circumstances, then the rank map of a locally of finite type projective module will be continuous. One of the main aims of the present paper is to realize this goal. In order to do this we need a suitable topology to replace instead of the ``discrete" topology of the natural numbers $\omega=\{0, 1, 2,...\}$. Finding such a topology requires to have a familiarity with the structure of the natural numbers. The set $\omega$ is an ordinal number and after some effort then we realized that the desired topology is the upper topology. First we studied the upper topology over an ordinal number and some new results are obtained, see Theorems \ref{Theorem II}, \ref{theorem 1} and \ref{th 190}. Then, in Theorem \ref{proposition 300}, it is shown that if  the ``discrete" topology of $\omega$ is replaced by the upper topology then the rank map remains continuous even if the projective module is infinitely generated. To prove its continuity the main ingredients which are used, in addition to the Kaplansky theorem, are some basic properties of the exterior powers of a module.\\

In this paper, we are also interested to investigate the continuity of the rank map of some certain modules. It should be noted that the projectivity of modules has very closed connection with the continuity of their rank maps. For instance, a f.g. flat module is projective if and only if its rank map is continuous with respect to the discrete topology. An analogue of this result is also proved for the upper topology, see Theorem \ref{theorem 721000}. Finally, Theorem \ref{proposition 1} is another interesting result of this paper which generalizes some major results in the literature on the projectivity of f.g. flat modules. \\

\section{Preliminaries}

We need the following material in the sequel. \\

Throughout the paper, all rings are commutative. \\

\begin{lemma}\label{lem 202} Let $(R,\mathfrak{m})$ be a local ring with the residue field $\kappa$ and let $M$ be a f.g. $R-$module. Then a finite subset $\{x_{1},...,x_{n}\}$ is a minimal generating set of $M$ if and only if $\{x_{i}+\mathfrak{m}M: 1\leq i\leq n\}$ is a $\kappa-$basis of $M/\mathfrak{m}M$.\\
\end{lemma}

{\bf Proof.} Well known. $\Box$ \\

\begin{corollary}\label{corollary 1} Any two minimal generating sets of a f.g. module over a local ring have the same number of elements. \\
\end{corollary}

{\bf Proof.} It is an immediate consequence of Lemma \ref{lem 202}. $\Box$ \\

The above Corollary does not hold in general. As a specific example, $\{1\}$ and $\{2,3\}$ are two minimal generating sets of $\mathbb{Z}$ as the module over itself with unequal number of elements.\\

Recall that for a given ring $R$ then the collection of subsets $D(f)\cap V(I)$ with $f\in R$ and $I$ runs through the set of f.g. ideals of $R$ is a basis for the opens of the patch (or, constructible) topology over $\Spec(R)$.  Moreover, the collection of subsets $V(I)$ where $I$ runs through the set of f.g. ideals of $R$ is a basis for the opens of the flat topology over $\Spec(R)$. The patch topology is clearly finer than the Zariski and flat topologies. The flat topology is the dual of the Zariski topology. For more details on the patch and flat topologies please consider \cite{Ebulfez}. \\

\begin{remark}\label{remark 1201} An $R-$module $M$ is called locally of finite type if $M_{\mathfrak{p}}$ as $R_{\mathfrak{p}}-$module is f.g. for all primes $\mathfrak{p}$. If $\mathfrak{p}$ is a prime ideal of $R$ then we define $\rn_{R_{\mathfrak{p}}}(M_{\mathfrak{p}})$ as the number of elements of a minimal generating set of the $R_{\mathfrak{p}}-$module $M_{\mathfrak{p}}$. By Corollary \ref{corollary 1}, it is well-defined. This leads us to a map from $\Spec(R)$ into the set of natural numbers $\omega=\{0,1,2,...\}$ given by $\mathfrak{p}\rightsquigarrow\rn_{R_{\mathfrak{p}}}(M_{\mathfrak{p}})$. It is called the rank map of $M$. Note that $\rn_{R_{\mathfrak{p}}}(M_{\mathfrak{p}})$, by Lemma \ref{lem 202}, is equal to the dimension of the $\kappa(\mathfrak{p})-$space $M\otimes_{R}\kappa(\mathfrak{p})$ where $\kappa(\mathfrak{p})$ is the residue field of $R$ at $\mathfrak{p}$. The rank map is said to be patch (resp. flat, Zariski) continuous if it is continuous whenever $\Spec(R)$ is equipped with the patch (resp. flat, Zariski) topology and $\omega$ with the discrete topology. It is easy to see that the rank map is locally constant if and only if it is patch continuous. The same statement is true for Zariski and also flat topologies. \\
\end{remark}

\begin{lemma}\label{lemma 2} Let $S$ and $T$ be two multiplicative subsets of $R$ such that $S\subseteq T$. Let $U$ be the image of $T$ under the canonical map $R\rightarrow S^{-1}R$. Then the ring $U^{-1}(S^{-1}R)$ is canonically isomorphic to $T^{-1}R$ and
for each $R-$module $M$, $U^{-1}(S^{-1}M)\simeq T^{-1}M$.\\
\end{lemma}

{\bf Proof.} It is well known. $\Box$ \\

Recall that a set $T$ is called a transitive set if $x\in T$ implies that $x\subset T$. A relation $<$ on a set $S$ is called linear (or, totally ordered) if it is transitive (i.e., if $x<y$ and $y<z$ then $x<z$) and satisfies the trichotomy law (i.e., for each pair $(x,y)$ of elements of $S$ then exactly one of the following conditions hold: $x<y$ or $x=y$ or $y<x$). We also say that $x\leq y$ if either $x=y$ or $x<y$. A linear relation $<$ over a set $S$ is said to be well-ordered if every non-empty subset of $S$ has the least element with respect to $<$. A set is called an ordinal (or, ordinal number) if it is a transitive set and with the membership relation $\in$ is well-ordered.\\

\begin{lemma}\label{lemma 565} Let $\alpha$ and $\beta$ be two ordinals. Then $\beta\in\alpha$ if and only if $\beta$ is a proper subset of $\alpha$. \\
\end{lemma}

{\bf Proof.} The implication ``$\Rightarrow$" is obvious. For the converse, see \cite[Theorem 7M, (d)]{Enderton}. $\Box$ \\

If $\alpha$ and $\beta$ are ordinals then we say that $\beta<\alpha$ if $\beta\in\alpha$ or equivalently $\beta\subset\alpha$.\\

\section{A remark on the Kaplansky theorem}

Let $M$ be a locally free $R-$module, i.e., $M_{\mathfrak{p}}$ is $R_{\mathfrak{p}}-$free for all $\mathfrak{p}\in\Spec(R)$.
We define $\rn_{R_{\mathfrak{p}}}(M_{\mathfrak{p}})$ as the cardinality of a $R_{\mathfrak{p}}-$basis of $M_{\mathfrak{p}}$. It is well-defined since over a commutative ring any two bases of a free module have the same cardinality. Note that projective modules and locally of finite type flat modules are typical examples of locally free modules.\\

\begin{lemma}\label{lemma 0108} Let $M$ be a locally free $R-$module and let $\mathfrak{p}\subseteq\mathfrak{q}$ be prime ideals of $R$. Then $\rn_{R_{\mathfrak{p}}}(M_{\mathfrak{p}})=
\rn_{R_{\mathfrak{q}}}(M_{\mathfrak{q}})$. \\
\end{lemma}

{\bf Proof.} By Lemma \ref{lemma 2}, $M_{\mathfrak{p}}\simeq(M_{\mathfrak{q}})_{\mathfrak{p}}$. $\Box$ \\

\begin{corollary}\label{lemma 0289} Let $M$ be a $R-$module. Then $\Supp M$ is stable under the specialization. If moreover $M$ is locally free then $\Supp M$ is stable under the generalization.\\
\end{corollary}

{\bf Proof.} Let $\mathfrak{p}\subseteq\mathfrak{q}$ be prime ideals of $R$.
If $\mathfrak{p}\in\Supp M$ then there is some $x\in M$ such that $\Ann_{R}(x)\subseteq\mathfrak{p}$. Thus $x/1$ is a non-zero element of $M_{\mathfrak{q}}$. Hence $\Supp M$ is stable under the specialization.
Now assume that $\mathfrak{q}\in\Supp M$. Then, by Lemma \ref{lemma 0108}, $M_{\mathfrak{p}}\neq0$. $\Box$ \\

\begin{lemma}\label{lemma 1} Let $M$ be a projective $R-$module and let $(I_{k})$ be a family of ideals of $R$. Then
$\bigcap\limits_{k}(I_{k}M)=
\big(\bigcap\limits_{k}I_{k}\big)M$.\\
\end{lemma}

{\bf Proof.} There is a free $R-$module $F$ such that $M$ is a direct summand of it. Thus there exists a $R-$submodule $N$ of $F$ such that $F=M+N$ and $M\cap N=0$. Let $\{x_{i}\}$ be a $R-$basis of $F$. Consider the bijective map $\psi:F\rightarrow\bigoplus\limits_{i}R$ given by $x\rightsquigarrow(r_{i})$ where $x=\sum\limits_{i}r_{i}x_{i}$. If $I$ is an ideal of $R$ then $\psi(IF)=\bigoplus\limits_{i}I$. Moreover $\bigcap\limits_{k}
(\bigoplus\limits_{i}I_{k})=
\bigoplus\limits_{i}
(\bigcap\limits_{k}I_{k})$. It follows that $\psi\big((\bigcap\limits_{k}I_{k})F\big)=
\psi\big(\bigcap\limits_{k}(I_{k}F)\big)$. Thus $(\bigcap\limits_{k}I_{k})F=\bigcap\limits_{k}(I_{k}F)$.
We also have $IF=IM+IN$ and $\bigcap\limits_{k}(I_{k}M+I_{k}N)=
\bigcap\limits_{k}(I_{k}M)+\bigcap\limits_{k}(I_{k}N)$.
Therefore
$(\bigcap\limits_{k}I_{k})M+
(\bigcap\limits_{k}I_{k})N=\bigcap\limits_{k}(I_{k}M)+
\bigcap\limits_{k}(I_{k}N)$. It follows that $(\bigcap\limits_{k}I_{k})M
=\bigcap\limits_{k}(I_{k}M)$. $\Box$ \\

\begin{theorem}\label{th 11} Let $M$ be either a projective $R-$module or a locally of finite type flat $R-$module and let $\mathfrak{p}$ be a prime ideal of $R$. Then $M_{\mathfrak{p}}=0$ if and only if $M=\mathfrak{p}M$.\\
\end{theorem}

{\bf Proof.} First assume that $M_{\mathfrak{p}}=0$. From the exact sequence: $$\xymatrix{0\ar[r]&R/\mathfrak{p}\ar[r]&\kappa(\mathfrak{p})}$$ we obtain the following exact sequence: $$\xymatrix{0\ar[r]&R/\mathfrak{p}\otimes_{R}M
\ar[r]&\kappa(\mathfrak{p})\otimes_{R}M.}$$
But $\kappa(\mathfrak{p})\otimes_{R}M\simeq R/\mathfrak{p}\otimes_{R}M_{\mathfrak{p}}=0$. It follows that $M=\mathfrak{p}M$. Conversely, assume that $M=\mathfrak{p}M$. We have $M_{\mathfrak{p}}$ as $R_{\mathfrak{p}}-$module is free because apply
the Kaplansky theorem \cite[Tag 0593]{Johan} whenever $M$ is projective and apply \cite[Tag 00NZ]{Johan} whenever $M$ is locally of finite type. Moreover $M_{\mathfrak{p}}\otimes_{R_{\mathfrak{p}}}\kappa(\mathfrak{p})
\simeq M_{\mathfrak{p}}\otimes_{R}R/\mathfrak{p}\simeq M/\mathfrak{p}M\otimes_{R}R_{\mathfrak{p}}=0$. It follows that $M_{\mathfrak{p}}=0$. $\Box$ \\

The following result is well known, see \cite[Lemma 6.2]{Lazard} or \cite[Lemma 1.1]{Vasconcelos}. We prove it by a new approach. \\

\begin{corollary}\label{th 22} The support of a projective module is Zariski $($and flat$)$ open.\\
\end{corollary}

{\bf Proof.} Let $M$ be a projective $R-$module, let $X=\Spec(R)\setminus\Supp(M)$ and let $I=\bigcap\limits_{\mathfrak{p}\in X}\mathfrak{p}$. By Theorem \ref{th 11} and Lemma \ref{lemma 1}, $IM=M$. Clearly $X\subseteq V(I)$. Conversely, assume that $I\subseteq\mathfrak{p}$. Then $M=IM\subseteq\mathfrak{p}M\subseteq M$. Thus $\mathfrak{p}M=M$ and so by Theorem \ref{th 11}, $\mathfrak{p}\in X$. Therefore $X=V(I)$. By using Corollary \ref{lemma 0289} and \cite[Theorem 3.11]{Ebulfez}, we conclude that $X$ is also flat closed. $\Box$ \\

\section{On the rank map of some certain modules}

\begin{lemma}\label{lemma 33} Let $M$ be a locally free $R-$module. Then for each natural number $n$, $\Supp\big(\Lambda^{n}(M)\big)=\{\mathfrak{p}\in\Spec(R): \rn_{R_{\mathfrak{p}}}(M_{\mathfrak{p}})\geq n\}$ where $\Lambda^{n}(M)$ is the $n$-th exterior power of $M$. \\
\end{lemma}

{\bf Proof.} By Theorem \ref{th 872}, $\Lambda^{n}(M)\otimes_{R}
R_{\mathfrak{p}}\simeq\Lambda_{R_{\mathfrak{p}}}^{n}(M_{\mathfrak{p}})$.
Therefore, by Corollary \ref{coro 5},  $\mathfrak{p}\in\Supp\Lambda^{n}(M)$ if and only if $\rn_{R_{\mathfrak{p}}}(M_{\mathfrak{p}})\geq n$.  $\Box$ \\

The following result is well known, e.g. see \cite[Corollary 2.2.2]{Weibel}. We provide a new proof for it. \\

\begin{proposition}\label{proposition 2} The rank map of a f.g. projective module is Zariski $($and flat$)$ continuous. \\
\end{proposition}

{\bf Proof.} Let $M$ be a f.g. projective module over a ring $R$. For each natural number $n$, then by Lemma \ref{lemma 33}, $\psi^{-1}(\{n\})=\Supp N\cap\big(\Spec(R)\setminus\Supp N'\big)$ where $\psi$ is the rank map of $M$, $N=\Lambda^{n}(M)$ and $N'=\Lambda^{n+1}(M)$. By Lemmas \ref{lemma 423} and \ref{lemma 457}, $N$ and $N'$ are f.g. projective $R-$modules. Thus there are idempotent elements $e,e'\in R$ such that $\Ann_{R}(N)=Re$ and $\Ann_{R}(N')=Re'$, because it is well known that the annihilator of every f.g. projective module is generated by an idempotent element. It follows that $\psi^{-1}(\{n\})=V(e)\cap V(1-e')=D\big(e'(1-e)\big)$.
Therefore $\psi^{-1}(\{n\})$ is both Zariski and flat open. $\Box$ \\

\begin{lemma}\label{lemma 566} Let $M$ be a locally free $R-$module which is also locally of finite type. Then the rank map of $M$ is patch continuous iff it is Zariski $($and flat$)$ continuous.\\
\end{lemma}

{\bf Proof.} For each natural number $n$, by Lemma \ref{lemma 0108}, $\psi^{-1}(\{n\})$ is stable under the generalization and specialization where $\psi$ is the rank map of $M$. Now if $\psi$ is patch continuous then, by \cite[Theorem 3.11]{Ebulfez}, it is both Zariski and flat continuous. The reverse is easy since the patch topology is finer than the Zariski topology. $\Box$ \\

\begin{theorem}\label{proposition 1} Let $R$ be a ring which has either a finitely many minimal primes or a finitely many maximal ideals. Then the rank map of a locally of finite type flat $R-$module is patch continuous. In particular, every f.g. flat $R-$module is $R-$projective.\\
\end{theorem}

{\bf Proof.} Let $M$ be a locally of finite type flat $R-$module, let $n$ be a natural number and let $E=\psi^{-1}(\{n\})$ where $\psi$ is the rank map of $M$. First assume that $\Min(R)=\{\mathfrak{p}_{1},...,\mathfrak{p}_{k}\}$. There exists some $s$ with $1\leq s\leq k$ such that $\mathfrak{p}_{s},\mathfrak{p}_{s+1},...,\mathfrak{p}_{k}\notin E$ but $\mathfrak{p}_{i}\in E$ for all $i<s$. By the proof of Lemma \ref{lemma 566},
$\Spec(R)\setminus E=\bigcup\limits_{i=s}^{k}V(\mathfrak{p}_{i})$. Thus $E$ is Zariski open in the case of finitely many minimal primes and so it is patch open. Now let $\MA(R)=\{\mathfrak{m}_{1},...,\mathfrak{m}_{d}\}$. Similarly, there exists some $\ell$ with $1\leq\ell\leq d$ such that $\mathfrak{m}_{\ell},\mathfrak{m}_{\ell+1},...,\mathfrak{m}_{d}\notin E$ but $\mathfrak{m}_{i}\in E$ for all $i<\ell$. Again by the proof of Lemma \ref{lemma 566}, $\Spec(R)\setminus E=\bigcup\limits_{i=k}^{d}\Lambda(\mathfrak{m}_{i})$ where $\Lambda(\mathfrak{m}_{i})=\{\mathfrak{p}\in\Spec(R) : \mathfrak{p}\subseteq\mathfrak{m}_{i}\}$. Therefore, by \cite[Corollary 3.6]{Ebulfez}, $E$ is a flat open in the case of finitely many maximal ideals and so it is patch open. The latter statement is an immediate consequence of \cite[Tags 00NZ, 00NX]{Johan} and Lemma \ref{lemma 566}.
$\Box$ \\

Theorem \ref{proposition 1} generalizes some major results in the literature on the projectivity of f.g. flat modules specially including \cite[Tag 00NZ]{Johan}, \cite[\S4E]{Lam}, \cite[Corollary 1.5]{Jondrup},
\cite[Fact 7.5]{Puninski-Rothmaler}. \\

\section{Upper topology and its applications}

Let $(P,<)$ be a poset. Then the collection of $d(x):=P\setminus\{y\in P: y\leq x\}$ with $x\in P$ is a sub-basis for the opens of the upper topology on $P$. Note that $\{y\in P: x<y\}\subseteq d(x)$. If the relation $<$ is linear (totally ordered) then the equality holds. The dual of the upper topology is called the lower topology. Thus the collection of $d'(x)=P\setminus\{y\in P: x\leq y\}$ with $x\in P$ is a sub-basis for the opens of the lower topology on $P$. Finally, the collection of $d(x)\cap d'(y)$ with $x,y\in P$ is a sub-basis for the opens of a topology on $P$. We call it the patch topology. \\

Let $R$ be a commutative ring and consider $\Spec(R)$ as a poset with respect to the strict inclusion. Then the lower topology over the poset $\Spec(R)$ is coarser than the Zariski topology. Similarly, the upper topology on the poset $\Spec(R)$ is coarser than the flat topology. \\

If $\mathbb{R}$ is the set of real numbers then the patch topology over the poset $(\mathbb{R}, <)$ coincides with the Euclidean topology. \\

Let $(P,<)$ be a poset and consider the upper topology over it. Then for each $x\in P$ we have $\overline{\{x\}}=\{y\in P: y\leq x\}$. If $S\subseteq P$ is a subset then the upper and subspace topologies over $S$ are the same. \\

Let $\alpha$ be an ordinal number. By the upper topology over $\alpha$ we mean the upper topology over the poset $(\alpha,\in)$. \\

\begin{theorem}\label{Theorem II} Let $\alpha$ be an ordinal and consider the upper topology over it. Then $\beta$ is a closed subset of $\alpha$ if and only if $\beta$ is an ordinal number with $\beta\leq\alpha$. \\
\end{theorem}

{\bf Proof.} If $F$ is a closed subset of $\alpha$ then to prove $F$ is an ordinal it will be enough to show that it is a transitive set. If $x\in F$ then
$\overline{\{x\}}\subseteq F$. We have $y\in\overline{\{x\}}$ if and only if either $y\in x$ or $y=x$. Hence, $F$ is a transitive set. Conversely, if $\beta$ is an ordinal with $\beta\leq\alpha$ then by Lemma \ref{lemma 565}, $\beta\subseteq\alpha$. If $x\in\beta$ then $\overline{\{x\}}\subseteq\beta$ because $\beta$ is a transitive set. Therefore $\beta$ is closed. $\Box$ \\

Recall that if $\alpha$ and $\beta$ are two ordinal numbers then $\alpha$ is called the successor of $\beta$ if $\alpha=\beta^{+}=\beta\cup\{\beta\}$.
An ordinal number is said to be a \emph{limit ordinal} if it is not the successor of an ordinal number. \\

\begin{theorem}\label{theorem 1} Let $\alpha$ be an ordinal and consider the upper topology over it. Then the following hold.\\
$\mathbf{(i)}$ If $\beta\in\alpha$ then $\overline{\{\beta\}}=\beta^{+}$.\\
$\mathbf{(ii)}$ $\alpha$ has a generic point if and only if it is not a limit ordinal. Moreover the generic point, if it exists, is unique.\\
$\mathbf{(iii)}$ If $\alpha\neq0$ then $\alpha$ is an irreducible space.\\
$\mathbf{(iv)}$ The closed subsets of $\alpha $ are stable under the arbitrary unions.\\
$\mathbf{(v)}$ Every open subset of $\alpha$ is quasi-compact. \\
$\mathbf{(vi)}$ $\alpha$ is a noetherian space.\\
\end{theorem}

{\bf Proof.} $\mathbf{(i)}$ and $\mathbf{(ii)}:$ It follows from the fact that an ordinal number is a transitive set. \\
$\mathbf{(iii)}:$ It follows from Theorem \ref{Theorem II} and the fact that any two ordinal numbers are comparable. \\
$\mathbf{(iv)}:$ It implies from Theorem \ref{Theorem II} and the fact that if $\{\beta_{i}\}$ is a family of ordinal numbers then $\bigcup\limits_{i}\beta_{i}$ is also an ordinal number, see e.g. \cite[Corollary 7N (d)]{Enderton}. \\
$\mathbf{(v)}:$ Let $U$ be an open subset of $\alpha$ and let $\{U_{i}\}_{i\in I}$ be an open covering of it. We may assume that $I\neq\emptyset$.
For each $i$ there exists an ordinal number $\beta_{i}$ with $\beta_{i}\leq\alpha$ such that $U_{i}^{c}=\beta_{i}$.
Let $\beta_{k}$ be the least element of the set $\{\beta_{i} : i\in I\}$ because it is well known that any non-empty set of ordinals has the least element.
It follows that $U=U_{k}$.\\
$\mathbf{(vi)}:$ It follows from (v). $\Box$ \\

\begin{remark}\label{remark 2} By \cite[Corollary 7N, (b) and (c)]{Enderton}, every natural number is an ordinal number where $0=\emptyset$, $1=0^{+}=\{0\}$, $2=1^{+}=\{0,1\}$ and so on. Moreover, by \cite[Corollary 7N, (a)]{Enderton} and \cite[Theorem 4G]{Enderton}, the set of natural numbers $\omega=\{0,1,2,...\}$ is also an ordinal number. In fact, $\omega$ is the first non-zero limit ordinal.\\
\end{remark}

The upper topology over an ordinal number $\alpha$ is discrete if and only if either $\alpha=0$ or $\alpha=1$. The upper topology over $2=\{0,1\}$ is just the Sierpi\'{n}ski topology.\\

\begin{theorem}\label{th 190} The upper topology over an ordinal number is spectral if and only if it is a natural number.\\
\end{theorem}

{\bf Proof.} First assume that the upper topology over an ordinal number $\alpha$ is spectral. If $\alpha\geq\omega$ then by Theorem \ref{Theorem II}, $\omega$ is a closed subset of $\alpha$ and so by Theorem \ref{theorem 1} (iii), it is irreducible. But $\omega$ does not have any generic point. Because, suppose there is some $m\in\omega$ such that $\omega=\overline{\{m\}}=m^{+}$. By \cite[Theorem 4I]{Enderton},
$m^{+}=m+1$ is a natural number, a contradiction. Hence $\alpha<\omega$, i.e., $\alpha$ is a natural number. Conversely, let $n$ be a natural number. To prove the assertion, by Theorem \ref{theorem 1}, it suffices to show that every closed and irreducible subset of $n$ has a generic point. We may assume that $n>0$. If $m$ is a closed and irreducible subset of $n$ then $m$ is a natural number with $0<m\leq n$. By \cite[Theorem 4C]{Enderton}, there exists a natural number $k$ such that $m=k^{+}$. Therefore $m=\overline{\{k\}}$. $\Box$ \\

As quoted in the Introduction, the rank map of a locally of finite type projective module is not necessarily continuous with respect to the discrete topology. But its continuity will be recovered if the discrete topology is replaced by the upper topology: \\

\begin{theorem}\label{proposition 300} Let $M$ be a locally of finite type flat $R-$module and consider the upper topology on the set of natural numbers. If $M$ is $R-$projective then the rank map of $M$ is Zariski $($and flat$)$ continuous.\\
\end{theorem}

{\bf Proof.} Let $n$ be a natural number. By \cite[Theorem 4C]{Enderton}, $n$ is either the empty set or
there exists a natural number $m$ such that $n=m+1$. Thus either $n=\emptyset$ or $n=\{0,1,2,...,m\}$. By applying this and Lemma \ref{lemma 33}, we obtain that $\psi^{-1}(n)=\psi^{-1}(\{0,1,2,...,m\})=
\Spec(R)\setminus\Supp\big(\Lambda^{n}(M)\big)$ where $\psi$ is the rank map of $M$. Therefore, by Lemma \ref{lemma 457} and Corollary \ref{th 22}, the assertion concludes. $\Box$ \\

\begin{theorem}\label{theorem 721000} Let $M$ be a f.g. flat module over a ring $R$. Then $M$ is $R-$projective if and only if the rank map of $M$ is patch continuous with respect to the upper topology.\\
\end{theorem}

{\bf Proof.} The implication ``$\Rightarrow$" implies from Proposition \ref{proposition 2} or Theorem \ref{proposition 300}. Conversely, let $n$ be a natural number. By Lemma \ref{lemma 33}, $\psi^{-1}(\{n\})=\Supp N\cap\big(\Spec(R)\setminus\Supp N'\big)$ and $\Supp N=\Spec(R)\setminus\psi^{-1}(n)$
where $\psi$ is the rank map of $M$, $N=\Lambda^{n}(M)$ and $N'=\Lambda^{n+1}(M)$. Thus, $\psi$ is patch continuous with respect to the discrete topology. Therefore, by \cite[Tags 00NZ, 00NX]{Johan} and Lemma \ref{lemma 566}, $M$ is $R-$projective. $\Box$ \\

In the proof of Theorem \ref{theorem 721000}, do not confuse $\psi^{-1}(\{n\})$ with $\psi^{-1}(n)$. \\

An $R-$module $M$ is said to be \emph{locally of countable rank} if for each prime ideal $\mathfrak{p}$ of $R$ then $M_{\mathfrak{p}}$ as $R_{\mathfrak{p}}-$module is countably generated (possibly infinite). Note that if $M$ is a projective $R-$module then $M=\bigoplus\limits_{i\in I}M_{i}$ is a direct sum of countably generated projective $R-$submodules, see \cite[Tag 058Y]{Johan}. If the index set $I$ is countable then $M$ is countably generated and so it is locally of countable rank.\\

\begin{corollary}\label{coro 6} Let $R$ be a ring such that $\Spec(R)$ is noetherian with respect to the flat topology. Let $M$ be a projective $R-$module which is locally of countable rank. Consider the upper topology over $\omega^{+}$ and the Zariski topology over $\Spec(R)$. Then the rank map of $M$ is continuous. \\
\end{corollary}

{\bf Proof.} We have $\psi^{-1}(\omega)=\bigcup\limits_{n\in\omega}\psi^{-1}(n)$ where $\psi$ is the rank map of $M$. By the proof of Theorem \ref{proposition 300}, $\psi^{-1}(n)=\Spec(R)\setminus\Supp\Lambda^{n}(M)$. Thus, by Corollary \ref{th 22}, it is Zariski closed. Therefore, by \cite[Theorem 4.2 ]{Ebulfez}, $\psi^{-1}(\omega)$ is also Zariski closed. $\Box$ \\

\begin{lemma}\label{lemma 400} Let $R$ be a ring such that $\Spec(R)$ is noetherian with respect to the Zariski topology. Then the flat opens of $\Spec(R)$ are stable under the arbitrary intersections.\\
\end{lemma}

{\bf Proof.} See \cite[Theorem 5.1]{Tarizadeh}. $\Box$ \\

As a dual of Corollary \ref{coro 6}, we have:\\

\begin{corollary} Let $R$ be a ring such that $\Spec(R)$ is noetherian with respect to the Zariski topology. Let $M$ be a projective $R-$module which is locally of countable rank. Consider the upper topology over $\omega^{+}$ and the flat topology over $\Spec(R)$. Then the rank map of $M$
is continuous. \\
\end{corollary}

{\bf Proof.} We have $\psi^{-1}(\omega)=\bigcup\limits_{n\in\omega}\psi^{-1}(n)$ where $\psi$ is the rank map of $M$. By Corollary \ref{th 22}, $\psi^{-1}(n)=\Spec(R)\setminus\Supp\Lambda^{n}(M)$ is a flat closed. Then apply Lemma \ref{lemma 400}. $\Box$ \\

All of the above results are symmetric by passing from the upper topology to the lower topology. This means that it does not matter with which topology (upper or lower) we work because we eventually get the same results.\\

\section{Appendix- exterior powers of a module}

Because of the lack of a good reference for our purposes and for the convenience of the reader, we present in this appendix some facts about the ``exterior powers of a module'' which are used in the body of the  paper. All of the presentation specially the proofs due to the author. The reader will also find this section a very useful introduction to the subject of ``exterior algebra". \\

Let $M, N$ be $R-$modules and let $n\geq2$. A function $f:M^{n}\rightarrow N$ is called \emph{alternative} if $f$ vanishes on each $n$-tuple with at least two distinct equal coordinates. The map $f$ is called \emph{skew-symmetric} if $f(x_{\sigma(1)},...,x_{\sigma(n)})=(\sgn\sigma)f(x_{1},...,x_{n})$ for all $\sigma\in S_{n}$. Here $S_{n}$ is the symmetric group of degree $n$. \\

\begin{lemma} Every alternative multi-linear map is skew-symmetric. \\
\end{lemma}

{\bf Proof.} Let $f:M^{n}\rightarrow N$ be a alternative multi-linear map and let $\sigma\in S_{n}$ where $n\geq2$. Each permutation can be written as a product of finitely many transpositions, see \cite[Lemma 6.7]{Isaacs}. Thus we may write $\sigma=\tau_{1}...\tau_{s}$ where $\tau_{k}$ is a transposition for all $k$. To prove the assertion we use an induction argument on $s$.
If $s=1$ then $\sigma=(i,j)$ with $i<j$. But $f(x_{1},...,x_{i}+x_{j},...,x_{j}+x_{i},...,x_{n})=0$. It follows that $f(x_{\sigma(1)},...,x_{\sigma(n)})=f(x_{1},...,x_{j},...,x_{i},...,x_{n})=
-f(x_{1},...,x_{i},...,x_{j},...,x_{n})$. Let $s>1$. We may write $\sigma=(i,j)\tau$. Using the induction hypothesis, then $f(x_{\sigma(1)},...,x_{\sigma(n)})=
f(x_{\tau(1)},...,x_{\tau(j)},...,x_{\tau(i)},...,x_{\tau(n)})=(\sgn\tau) f(x_{1},...,x_{j},...,x_{i},...,x_{n})=
(\sgn\sigma)f(x_{1},...,x_{i},...,x_{j},...,x_{n})$. $\Box$ \\

\begin{lemma}\label{lemma 3} A multi-linear map $f:M^{n}\rightarrow N$ is alternative if and only if the map $f$ vanishes on each $n$-tuple with a pair of adjacent equal coordinates.\\
\end{lemma}

{\bf Proof.} Suppose for a given $n$-tuple $(x_{1},...,x_{n})\in M^{n}$, $x_{i}=x_{j}$ with $i<j$. If $j=i+1$ then there is nothing to prove. Let $j>i+1$. By the hypothesis, $f(x_{1},...,x_{i},...,x_{j-1}+x_{j}, x_{j}+x_{j-1},...,x_{n})=0$.
By the induction hypothesis, $f(x_{1},...,x_{i},...,x_{j},x_{j-1},...,x_{n})=0$. It follows that $f(x_{1},...,x_{i},...,x_{j-1},x_{j},...,x_{n})=0$. $\Box$ \\

Let $J_{n}$ be the $R-$submodule of $M^{\otimes n}$ generated by the collection of pure tensors of the form $x_{1}\otimes...\otimes x_{n}$ with $x_{i}=x_{j}$ for some $i\neq j$. The quotient module  $\Lambda^{n}(M):=M^{\otimes n}/J_{n}$ is called the $n$-th exterior power of $M$. Write $\Lambda^{0}(M)=R$ and $\Lambda^{1}(M)=M$.
The canonical multi-linear map $\eta:M^{n}\rightarrow\Lambda^{n}(M)$ given by $(x_{1},...,x_{n})\rightsquigarrow x_{1}\wedge...\wedge x_{n}:=x_{1}\otimes...\otimes x_{n}+J_{n}$ is clearly alternative and the module $\Lambda^{n}(M)$ together with the map $\eta$ satisfy in the following universal property: for each alternative multi-linear map $\psi:M^{n}\rightarrow N$ then there is a unique $R-$linear map $\phi:\Lambda^{n}(M)\rightarrow N$ such that $\psi=\phi\circ\eta$. \\

Specially, consider the map $M^{n}\rightarrow M^{\otimes n}$ given by $(x_{1},...,x_{n})\rightsquigarrow\sum\limits_{\sigma\in S_{n}}
(\sgn\sigma)x_{\sigma(1)}\otimes...\otimes x_{\sigma(n)}$. This map
is multi-linear since multi-linear maps are stable under the finite sums.
It is also alternative. Because suppose $x_{i}=x_{i+1}$ for some $i$. Let $\sigma\in S_{n}$. We have $\sigma(j)=i$ and $\sigma(k)=i+1$ for some $j,k$.
Consider the permutation $\sigma^{\ast}=(j,k)\sigma$. Then clearly $\sgn\sigma^{\ast}=-\sgn\sigma$, $x_{\sigma(1)}\otimes...\otimes x_{\sigma(n)}=x_{\sigma^{\ast}(1)}\otimes...\otimes x_{\sigma^{\ast}(n)}$ and the assignment $\sigma\rightsquigarrow\sigma^{\ast}$ is an injective map from $A_{n}$, the set of even permutations of degree $n$, onto the set of odd permutations of the same degree. It follows that $\sum\limits_{\sigma\in S_{n}}
(\sgn\sigma)x_{\sigma(1)}\otimes...\otimes x_{\sigma(n)}=\sum\limits_{\sigma\in A_{n}}
(x_{\sigma(1)}\otimes...\otimes x_{\sigma(n)}-
x_{\sigma^{\ast}(1)}\otimes...\otimes x_{\sigma^{\ast}(n)})=0$. Therefore, by the universal property of the exterior powers, there is a unique morphism of $R-$modules $\delta:\Lambda^{n}(M)\rightarrow M^{\otimes n}$ which maps each pure wedge $x_{1}\wedge...\wedge x_{n}$ into $\sum\limits_{\sigma\in S_{n}}
(\sgn\sigma)x_{\sigma(1)}\otimes...\otimes x_{\sigma(n)}$. The map $\delta$ is injective whenever $M$ is a free $R-$module. Because if $\{x_{\alpha} : \alpha\in I\}$ is a basis of $M$ then the collection of pure tensors $x_{i_{1}}\otimes...\otimes x_{i_{n}}$ with $\{i_{1},...,i_{n}\}\subseteq I$ is a basis of $M^{\otimes n}$. In particular, we have shown that:\\

\begin{lemma}\label{lemma 423} Let $M$ be a (resp. free) $R-$module and let $\{x_{\alpha} : \alpha\in I\}$ be a generating set (resp. basis) of $M$. Consider a well-ordering relation $<$ on the index set $I$. Then the collection of pure wedges of the form $x_{i_{1}}\wedge...\wedge x_{i_{n}}$ where $\{i_{1},...,i_{n}\}\subseteq I$ with $i_{1}<...<i_{n}$ is a generating set (resp. basis) of the $R-$module $\Lambda^{n}(M)$. $\Box$ \\
\end{lemma}

\begin{corollary}\label{coro 5} If $M$ is a free $R-$module of rank $n$ then $\Lambda^{i}(M)$ is a free $R-$module of rank $\binom{n}{i}$ for $0\leq i\leq n$ and $\Lambda^{i}(M)=0$ for all $i>n$. $\Box$ \\
\end{corollary}

Let $\phi:M\rightarrow N$ be a $R-$linear map. By the universal property of the exterior powers, then there is a unique $R-$linear map $\Lambda^{n}(\phi):\Lambda^{n}(M)\rightarrow\Lambda^{n}(N)$ which maps each pure wedge $x_{1}\wedge...\wedge x_{n}$ into $\phi(x_{1})\wedge...\wedge\phi(x_{n})$. In fact, $\Lambda^{n}(-)$ is an additive functor from the category of $R-$modules into itself.\\

\begin{lemma}\label{lemma 4588} Let $(M_{i}, \phi_{i,j})$ be an inductive (direct) system of $R-$modules over the directed set $(I,<)$. Then for each fixed $n\geq0$, $\big(\Lambda^{n}(M_{i}), \Lambda^{n}(\phi_{i,j})\big)$ is an inductive system of $R-$modules over the same directed set and we have the following canonical isomorphism of $R-$modules:  $$\lim\limits_{\overrightarrow{i\in I}}\Lambda^{n}(M_{i})\simeq\Lambda^{n}(\lim
\limits_{\overrightarrow{i\in I}}M_{i}).$$ \\
\end{lemma}

{\bf Proof.} Easy. $\Box$ \\

\begin{lemma}\label{lemma 457} If $M$ is $R-$projective (resp. $R-$flat) then $\Lambda^{n}(M)$ is as well.\\
\end{lemma}

{\bf Proof.} First assume that $M$ is $R-$projective. Consider the exact sequence  $\xymatrix{0\ar[r]&N\ar[r]^{\phi}&F\ar[r]^{\psi}&M\ar[r]&0}$ where $F$ is a free $R-$module. It is split since $M$ is projective. In an abelian category, every exact and split sequence is left exact and split by an additive functor. Therefore the following sequence is exact and split $$\xymatrix{0\ar[r]&\Lambda^{n}(N)
\ar[r]^{\Lambda^{n}(\phi)}&\Lambda^{n}(F)
\ar[r]^{\Lambda^{n}(\psi)}&\Lambda^{n}(M)\ar[r]&0.}$$ It follows that $\Lambda^{n}(F)\simeq\Lambda^{n}(M)\oplus\Lambda^{n}(N)$.
By Lemma \ref{lemma 423}, $\Lambda^{n}(F)$ is a free $R-$module. Thus $\Lambda^{n}(M)$ is $R-$projective. If $M$ is a flat $R-$module then the flatness of $\Lambda^{n}(M)$ is an immediate consequence of Corollary \ref{coro 5}, Lemma \ref{lemma 4588} and \cite[Tag 058G]{Johan}. $\Box$ \\

\begin{theorem}\label{th 872} Let $R\rightarrow S$ be a morphism of rings and let $M$ be an $R-$module. Then $\Lambda^{n}(M)\otimes_{R}S$ as $S-$module is canonically isomorphic to $\Lambda_{S}^{n}(M\otimes_{R}S)$. \\
\end{theorem}

{\bf Proof.} Let $N=M\otimes_{R}S$. Consider the canonical map $M^{n}\rightarrow \Lambda_{S}^{n}(N)$ which maps each $n$-tuple $(x_{1},...,x_{n})\in M^{n}$ into $(x_{1}\otimes1)\wedge...\wedge(x_{n}\otimes1)$. Clearly it is $R-$multilinear and alternative. Therefore there is a (unique) morphism of $R-$modules $\Lambda^{n}(M)\rightarrow\Lambda_{S}^{n}(N)$ which maps each pure wedge $x_{1}\wedge...\wedge x_{n}$ of $\Lambda^{n}(M)$ into $(x_{1}\otimes1)\wedge...\wedge(x_{n}\otimes1)$. Then we obtain a morphism of $S-$modules $\phi:\Lambda^{n}(M)\otimes_{R}S\rightarrow\Lambda_{S}^{n}(N)$ which maps each pure tensor $(x_{1}\wedge...\wedge x_{n})\otimes s$ into
$(x_{1}\otimes s)\wedge(x_{2}\otimes1)\wedge...\wedge(x_{n}\otimes1)$.
In order to find the inverse of $\phi$ we act as follows. Consider the canonical map $f:N^{n}\rightarrow\Lambda^{n}(M)\otimes_{R}S$ which maps each $n$-tuple of pure tensors $(x_{1}\otimes s_{1},...,x_{n}\otimes s_{n})$ into $(x_{1}\wedge...\wedge x_{n})\otimes s_{1}...s_{n}$. Clearly it is $S-$multilinear.
We claim that it is also an alternative map. Because by Lemma \ref{lemma 3}, it suffices to prove that $f(z_{1},...,z_{n})=0$ whenever  $z_{k}=z_{k+1}=\sum\limits_{i=1}^{d}x_{i}\otimes s_{i}$ for some $k$. But $f(z_{1},...,z_{n})$ is a finite sum of elements of the form $\sum\limits_{1\leq i,j\leq d}(x'_{1}\wedge...\wedge x'_{k-1}\wedge x_{i}\wedge x_{j}\wedge x'_{k+2}\wedge...\wedge x'_{n})\otimes s'_{1}...s'_{k-1}s_{i}s_{j}s'_{k+2}...s'_{n}$.
We show that each of them is zero. In fact, it suffices to verify that $\sum\limits_{1\leq i,j\leq d}(x_{i}\wedge x_{j})\otimes s_{i}s_{j}=0$.
We act by induction on $d$. If $d\geq2$ then we may write
$\sum\limits_{1\leq i,j\leq d}(x_{i}\wedge x_{j})\otimes s_{i}s_{j}=
\sum\limits_{j=1}^{d-1}\Big(\sum\limits_{i=1}^{d-1}(x_{i}\wedge x_{j})\otimes s_{i}s_{j}+(x_{d}\wedge x_{j})\otimes s_{d}s_{j}\Big)+\sum\limits_{i=1}^{d}(x_{i}\wedge x_{d})\otimes s_{i}s_{d}=
\sum\limits_{j=1}^{d-1}(x_{d}\wedge x_{j})\otimes s_{d}s_{j}-\sum\limits_{i=1}^{d-1}(x_{d}\wedge x_{i})\otimes s_{d}s_{i}=0$. This establishes the claim. Therefore there is a (unique) morphism of $S-$modules $\psi:\Lambda^{n}(N)\rightarrow\Lambda^{n}(M)\otimes_{R}S$ which maps each pure wedge of the form $(x_{1}\otimes s_{1})\wedge...\wedge(x_{n}\otimes s_{n})$ into $(x_{1}\wedge...\wedge x_{n})\otimes s_{1}...s_{n}$. Then it is easy to see that $\phi\circ\psi$ and $\psi\circ\phi$ are the identity. $\Box$ \\

\end{document}